\documentclass[11pt,oneside]{amsart}
\usepackage{amsmath,ifthen, amsfonts, amssymb, graphicx, srcltx, mathrsfs, xfrac}

\newcommand{\showcomments}{yes}

\newsavebox{\commentbox}
{\ifthenelse{\equal{\showcomments}{yes}}%
{\footnotemark
    \begin{lrbox}{\commentbox}
    \begin{minipage}[t]{1.25in}\raggedright\sffamily\tiny
    \footnotemark[\arabic{footnote}]}
{\begin{lrbox}{\commentbox}}}%
{\ifthenelse{\equal{\showcomments}{yes}}%
{\end{minipage}\end{lrbox}\marginpar{\usebox{\commentbox}}}
{\end{lrbox}}}

\newcommand{\diam}{\mathrm{diam}}
\newcommand{\Area}{\mathrm{Area}}

\newcommand{\diff}{\mathrm{d}}

\newtheorem{thm}{Theorem}[section]
\newtheorem{lem}[thm]{Lemma}
\newtheorem{cor}[thm]{Corollary}

\newtheorem{prop}[thm]{Proposition}

\theoremstyle{definition}
\newtheorem{defn}[thm]{Definition}
\newtheorem{rem}[thm]{Remark}

\newcommand{\field}[1]{\mathbb{#1}}
\newcommand{\integers}{\ensuremath{\field{Z}}}
\newcommand{\naturals}{\ensuremath{\field{N}}}
\newcommand{\reals}{\ensuremath{\field{R}}}

      \title{A new length estimate for curve shortening flow and low regularity initial data}

\author{Joseph Lauer}

\begin{document}

\date{November 20, 2011}

\begin{abstract} 
In this paper we introduce a geometric quantity, the $r$-multiplicity, that controls the length of a smooth curve as it evolves by curve shortening flow.  The length estimates we obtain are used to prove results about the level set flow in the plane. If  $K$ is locally-connected, connected and compact, then the level set flow of $K$ either vanishes instantly, fattens instantly or instantly becomes a smooth closed curve.  If the compact set in question is a Jordan curve $J$, then the proof proceeds by using the $r$-multiplicity to show that if $\gamma_n$ is a sequence of smooth curves converging uniformly to $J$, then the lengths $\mathscr{L}({\gamma_n}_t)$, where ${\gamma_n}_t$
denotes the result of applying curve shortening flow to $\gamma_n$ for time t, are uniformly bounded for each $t>0$.  Once the level set flow has been shown to be smooth we prove that the Cauchy problem for curve shortening flow has a unique solution if the initial data is a finite length Jordan curve.    
\end{abstract}

\maketitle

\section{Introduction}

A smooth map $\gamma:\mathbb{M}\times(t_1,t_2)\to\reals^2$, where $\mathbb{M}$ is either $\reals$ or $\mathbb{S}^1\cong\frac{\reals}{2\pi\integers}$, is a solution to curve shortening flow if 
$$
\frac{\partial\gamma}{\partial t}=\kappa\vec{n},
$$
where $\kappa\vec{n}$ is the well-defined curvature vector.  Geometrically, a solution produces a 1-parameter family of evolving curves which we denote by $\gamma_t$. 

Short-time existence for smooth initial data in the case $\mathbb{M}=\mathbb{S}^1$ was proved by Gage and Hamilton~\cite{GH86}.  In that paper, it was also shown that any convex curve shrinks to a point, becoming asymptotically round, and Grayson~\cite{G87} later proved that an arbitrary smooth closed curve becomes convex.

Curve shortening flow (CSF) bears a strong connection to the heat equation, and although it is nonlinear it exhibits smoothing properties common in parabolic PDE's.   In particular, like the heat equation, CSF is able to smooth non-smooth initial data:  The first results in this direction were due to Ecker and Huisken~\cite{EH89} who proved that if $\gamma$ is a smooth entire Lipschitz graph, then for each $t>0$, the curvature of $\gamma_t$, and all its derivatives, are bounded in terms of $t$ and the Lipschitz constant of the initial data.  By approximation this leads to the result ~\cite{EH89} that the curve shortening flow has a smooth solution whenever the initial data is an entire Lipschitz graph.  

In~\cite{EH91} the same authors proved that the same conclusion is true for entire graphs which are merely locally-Lipschitz, and gave a so-called uniformly locally-Lipschitz condition which provides an existence theorem that can be applied to closed curves.  We note that the results in ~\cite{EH89} and~\cite{EH91} apply more generally to mean curvature flow in $\reals^n$.
 
One of the major obstacles to extending the results above in the case of closed curves is that analytic estimates are difficult to control when the lengths of any approximating sequence are unbounded.  In this paper we establish a length estimate that is strong enough to allow for applications to the evolution of nonsmooth sets.  By approximating a locally-connected, compact set $K$ by a sequence of smooth curves whose $r$-multiplicity stays bounded, we are able to prove that the level set flow of $K$, denoted by $K_t$, instantly becomes smooth.  The level set flow is a weak notion of CSF (or more generally, of mean curvature flow) that allows one to evolve an arbitrary compact set in a way that agrees with CSF when the initial data is a smooth embedded curve.  For the level set flow we obtain:

\begin{thm} \label{main} Let $K\subset\reals^2$ be locally-connected, connected and
compact. Then there exists $T>0$ such that exactly one of the following holds for $0<t<T$:
\newline\indent (1) $K_t=\emptyset$.
\newline\indent (2) $K_t$ is a smooth closed curve.
\newline\indent (3) The interior of $K_t$ is nonempty. 
\end{thm}

Moreover, which of the three categories a particular set $K$ falls into depends only on the number of components of $\reals^2\setminus K$, and the Lebesgue measure of $K$, denoted by $m(K)$.  Even in the case when $K_t$ fattens the level set flow is well understood since we show that $\partial K_t$ consists of finitely many disjoint smooth closed curves for each $t>0$.  For Jordan curves the situation depends on the Lebesgue measure of the curve.  

\begin{thm}\label{jordan_curve}  Let $J\subset\reals^2$ be a Jordan curve.  Then there exists $T>0$ such that for $0<t<T$
\newline\indent (1) if $m(J)=0$, then $J_t$ is a smooth closed curve, and
\newline\indent (2) if $m(J)>0$, then $J_t$ is an annular region with smooth boundary.
\end{thm}

The fact that the level set flow is smooth leads to the following existence theorem for CSF, which is stated precisely in $\S$\ref{cauchy_section}.  We note that this case includes curves which are not locally graphs.

\begin{thm} The Cauchy problem for curve shortening flow has a unique solution if the initial data is a finite length Jordan curve.  
\end{thm}

\subsection{Length estimates using $r$-multiplicity} The $r$-multiplicity, which we denote by $M_r(u)$, is defined for all continuous maps $u:\mathbb{S}^1\to\reals^2$ and $r>0$.  See Section~\ref{rmult_section} for the defintion.  For our purposes the usefulness of the $r$-multiplicity stems from the fact that it is comparable to the quantity $\widetilde{M}_r(u)$ defined below.

Let $\mathcal{N}_r(\ell)$ denote the open $r$-neighbourhood of~$\ell$.

\smallskip

\begin{figure}
\centering
\scalebox{0.8}{\includegraphics{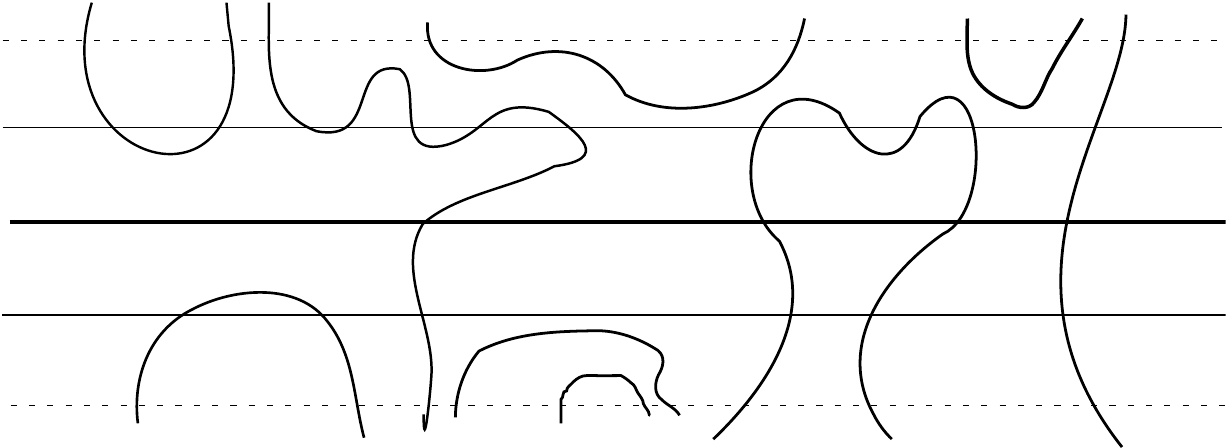}}
\caption{The central horizontal line is $\ell$, and the other lines are the boundaries of the strips of radius $r$ and $2r$ around $\ell$.  In this example $\gamma\cap\mathcal{N}_{2r}(\ell)$ contains 9 components and $M_{r,\ell}(\gamma)=5$.}  
\label{mult_def_fig}
\end{figure}

\begin{defn} ~\label{intro_mult_def} 
Fix $r>0$ and let $\ell$ be a line in $\reals^2$.  Given a continuous map $u:\mathbb{S}^1\to\reals^2$,  let $\{I_i\}_{i\in\Lambda}$ be the connected components of $u^{-1}(u(\mathbb{S}^1)\cap\mathcal{N}_{2r}(\ell))$ such that $u(I_i)\cap\overline{\mathcal{N}_r(\ell)}\neq\emptyset$.  Then define $\widetilde{M}_{r,\ell}(u)=|\Lambda|$, and
$$
\widetilde{M}_r(u)=\sup_{\mathrm{lines}\:\:\ell}\{\widetilde{M}_{r,\ell}(u)\}.
$$
\end{defn}

See Figure~\ref{mult_def_fig} for an example of the calculation of $\widetilde{M}_{r,\ell}(u)$ in the case that $u$ is an embedding. 

We now explain in what sense the $r$-multiplicity controls the length of a curve as it evolves.  Given $t>0$ and a smooth closed curve $\gamma$, there is a particular scale $r$, so that $\mathscr{L}(\gamma_t)$ is controlled by the $r$-multiplicity of $\gamma=\gamma_0$ at that scale.  Our main estimate is:

\begin{thm} \label{main_estimate_intro}
For all $t,d>0$ there exists constants $r(t,d)$ and
$C=C(t,d)$ with the following property:  Let
$\gamma$ be a smooth embedded closed curve with $\diam (\gamma)<d$.
Then
$$
\mathscr{L}(\gamma_t)<CM_r(\gamma).
$$
\end{thm}

Although this statement includes only embedded curves, our methods apply to immersed curves, and the most general statement can be found in Theorem~\ref{main_estimate}.

In Section~\ref{rmult_section} we observe that the $r$-multiplicity is uniformly bounded under uniform convergence and so we obtain the following, which is the basis for our results concerning the level set flow of Jordan curves.

\begin{cor} \label{convergence_estimate_intro} For each continuous map $u:\mathbb{S}^1\to\reals^2$ and $t>0$, there is
a  constant $\widetilde{C}=\widetilde{C}(u,t)$ such that if $\gamma_n$ is any sequence of smooth embedded 
curves which converge uniformly to $u$, then
$$
\mathscr{L}({\gamma_n}_t)<\widetilde{C}
$$
for sufficiently large $n$.
\end{cor}

\begin{figure}
\centering
\scalebox{0.75}{\includegraphics{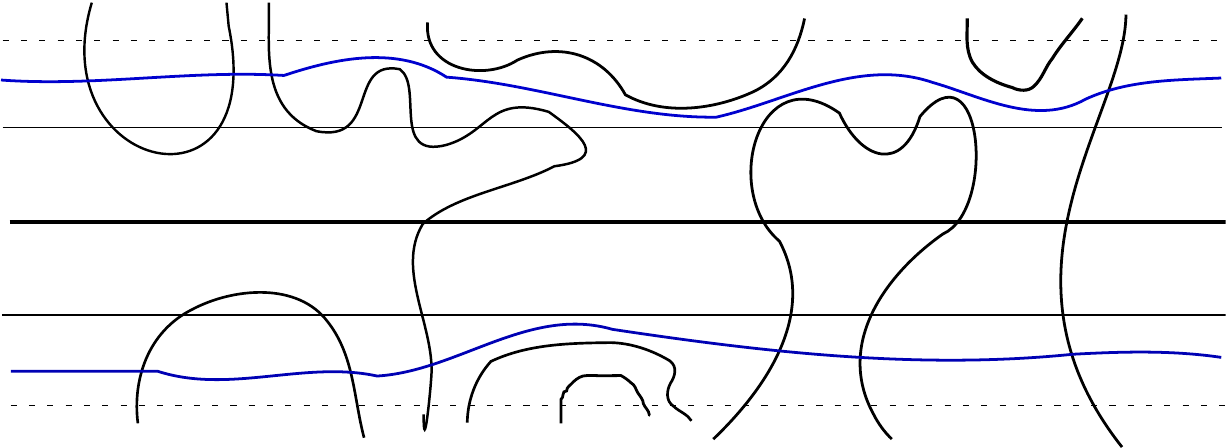}}
\caption{The outermost leaves of the foliation $\mathcal{F}$ have been added to Figure~\ref{mult_def_fig}.  Each leaf intersects each component of $\gamma\cap\mathcal{N}_{2r}(\ell)$ that counts towards $\widetilde{M}_{r,\ell}(\gamma)$ at most twice and is disjoint from the other components.}
\label{foliate_fig}
\end{figure}

If such a sequence exists then $u$ is called a near-embedding.  This class includes certain area-filling curves, and it is natural to ask to what extent a sequence of smooth approximations can be used to define an evolution of a near-embedding.

\subsection{Outline of the proof of Theorem~\ref{main_estimate_intro}}\label{intro_outline} For a moment let us assume that we have the function~$r(t,d)$ in Theorem~\ref{main_estimate_intro}.  We will explain shortly how this function is determined.  Given a smooth closed curve $\gamma$ and $t>0$, we proceed by bounding $\mathscr{L}(\gamma_t\cap B_r(x))$, in terms of the $r$-multiplicity, at each $x\in\reals^2$.  The global bound follows immediately since the diameter of $\gamma$ is controlled. 

Let~$\ell$~be a line through $x\in\reals^2$.   In Lemma~\ref{foliate1}, we construct a smooth foliation, $\mathcal{F}$, of a region  containing the strip $\mathcal{N}_r(\ell)$.  $\mathcal{F}$ is constructed  so that~(1)~each leaf is parallel to~$\ell$~outside a small neighbourhood of  $\gamma$, and~(2)~each leaf intersects $\gamma$ at most $2\widetilde{M}_r(\gamma)$ times.  See Figure~\ref{foliate_fig}.  By evolving each leaf in $\mathcal{F}$ by CSF, we obtain a foliation~$\mathcal{F}_t$ for each $t>0$.  In Section \ref{straight_section}, we prove the following straightening lemma which is applied to each leaf of $\mathcal{F}$.

\begin{lem} \label{straightening_intro} Given $r,l,\alpha>0$, there exists
$T=T(r,l,\alpha)>0$ with the following property: If $\gamma$ is a smooth
curve whose image is equal to the x-axis outside
$[0,l]\times[-r,r]$, then $\gamma_t$ is an $\alpha$-Lipschitz graph for
all $t\geq T$.

Moreover, with $l$ and $\alpha$ fixed we have $T\to 0$ as $r\to 0$.
\end{lem}

\begin{figure}
\centering
\scalebox{0.8}{\includegraphics{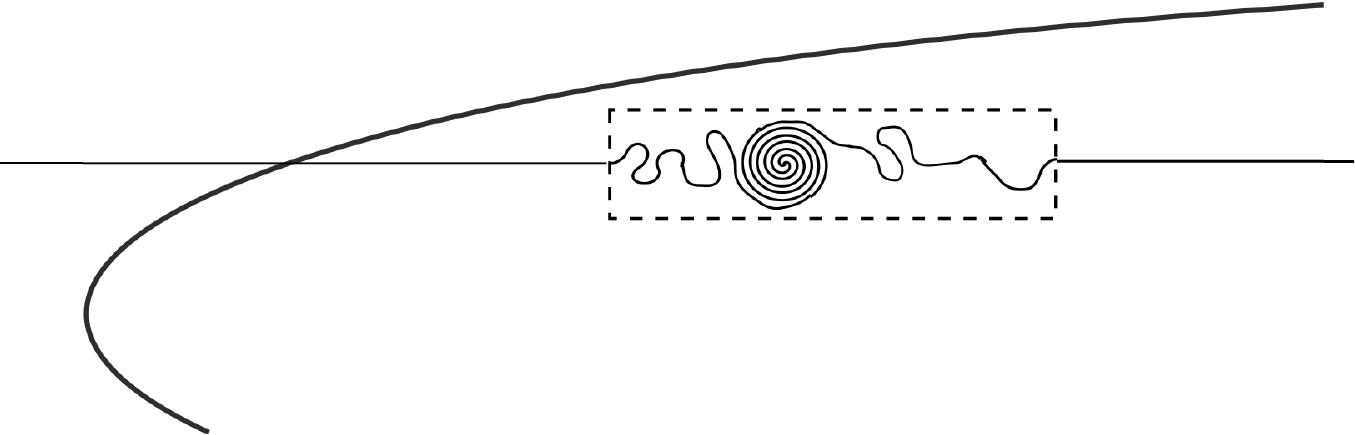}}
\caption{The proof of Lemma~\ref{straightening_intro} uses grim reapers that have bounded slope at all points that pass through the rectangle $[0,l]\times[-r,r]$.  Since the number of intersections does not increase, the slope of $\gamma_t$ will be bounded as soon as a single grim reaper has passed completely through the rectangle.}
\label{Reaper2}
\end{figure}

The proof of Lemma~\ref{straightening_intro} uses a family of grim reapers that intersect $\gamma=~\gamma_0$ exactly once, and travel in the direction of the $x$-axis with speed proportional to $r^{-1}$.  Moreover, we choose the family of grim reapers so that the slope of the tangent is bounded by $\alpha$ at any point that passes through $[0,l]\times[-r,r]$.  See Figure~\ref{Reaper2}.  After some explicitly given time $T$, $[0,l]\times[-r,r]$ is foliated by segments of grim reapers that intersect $\gamma_T$ exactly once, and this allows us to bound the slope of the tangent of $\gamma_T$.  

The final statement of Lemma~\ref{straightening_intro} allows us to define the function $r(t,d)$:  That is, given a time $t_0>0$ that we would like to make conclusions about, we can choose $r>0$ small enough so that $t_0=T(2r,l,\alpha)$, where $0<\alpha\ll 1$ is a small constant and $l=\diam(\gamma)+2$.  With $r$ chosen in this way Lemma~\ref{straightening_intro} implies that $\mathcal{F}_{t_0}$ consists of $\alpha$-Lipschitz graphs over~$\ell$.   

We then repeat this procedure with the line through $x$ perpendicular to~$\ell$, obtaining a second foliation $\mathcal{F}'_{t_0}$.  Both $\mathcal{F}_{t_0}$ and $\mathcal{F}'_{t_0}$ cover $B_r(x)$, and intersect transversely there.  Moreover, $\gamma_{t_0}$ intersects each leaf in $\mathcal{F}_{t_0}$ and $\mathcal{F}'_{t_0}$ at most $2\widetilde{M}_r(\gamma)$ times since the number of intersections does not increase along the flow.  In Section~\ref{holonomy_section} we show that the restriction of $\mathcal{F}_{t_0}$ and $\mathcal{F}'_{t_0}$ to $B_r(x)$ is bilipschitz equivalent to the standard grid, and that the bilipschitz constant depends only on $t_0$~and~$r$.  This follows from the fact that the foliations have been evolving by CSF and the linearization of CSF plays a key role.  Thus we obtain estimates for $\mathscr{L}(\gamma_{t_0}\cap B_r(x))$ in terms of the length bounds for a smooth curve that intersects each horizontal and vertical line a controlled number of times.

\subsection{Outline of sections}  

In $\S$\ref{rmult_section},  the basic properties of $r$-multiplicity are proved.  In $\S$\ref{straight_section} and $\S$\ref{holonomy_section} we explain how the foliations described above can be sufficiently smoothed.  In $\S$\ref{straight_section} we show that each leaf straightens individually, while in $\S$\ref{holonomy_section} we show that the holonomy maps of the evolving foliaton are controlled at all positive times.  In $\S$~\ref{construct_section} the initial foliation is constructed and Theorem~\ref{main_estimate_intro} is proved in $\S$\ref{length_estimates_section}.  In $\S$\ref{lc_sets_section} we prove the existence of bounded $r$-multiplicity approximations to locally-connected sets which are not curves.  In $\S$\ref{LSF_section}  we discuss the level set flow in the plane, and give a characterization of it in terms of sequences of smooth curves, and in $\S$\ref{smoothness_section} we apply the length estimates of $\S$\ref{length_estimates_section} to show that the level set flow is smooth.  

In the last two sections we restrict to the case of Jordan curves.  In $\S$\ref{backwards_convergence_section} it is shown that the level set flow converges onto the original curve without overlaps, and in $\S$\ref{cauchy_section} we show that the Cauchy problem for CSF has a solution for finite length initial data.

\subsection{Acknowledgements} The author wishes to thank Bruce Kleiner for his guidance and direction during the project, and for suggesting the problem.


\section{$r$-multiplicity}~\label{rmult_section}

In this section we introduce the $r$-multiplicity and collect its most basic properties.

Recall that $\mathcal{N}_r(\ell)$ denotes the open $r$-neighbourhood of $\ell$.

\begin{defn} [$r$-multiplicity] ~\label{rmult_new} Fix $r>0$ and let $\ell$ be a line in $\reals^2$.  Given a continuous map $u:\mathbb{S}^1\to\reals^2$, let $\{I_i\}_{\in\Delta}$
be the connected components of  $u^{-1}(\mathcal{N}_{r/2}(\ell)\cap u(\mathbb{S}^1))$ such that $\overline{u(I_i)}$ intersects both components of $\partial\mathcal{N}_{r/2}(\ell)$, and define $M_{r,\ell}(u)=|\Delta|$. 

We define the {\it $r$-multiplicity} of $u$ by
$$
M_r(u)=\sup_{\mathrm{lines}\:\:\ell}\{M_{r,\ell}(u)\}.
$$ 
\end{defn}

Intuitively, $M_r(u)$ is the maximum number of times that $u$ crosses an infinite strip of height~$r$.  It is comparable to the quantity $\widetilde{M}_r(u)$ defined in Definition~\ref{intro_mult_def}, but is more natural, and easier to use in proofs.  It is easy to check that 
$$
\frac{1}{2}M_r(u)\leq\widetilde{M}_r(u)\leq M_r(u).
$$

The $r$-multiplicity may be thought of as a coarse version of intersection number.  However, unlike intersections it is always finite:

\begin{prop}\label{multiplicity}
$M_r(u)<\infty$.
\end{prop}

\proof This follows from the uniform continuity of the function $u:\mathbb{S}^1\to\reals^2$.  Indeed, given $r>0$ there exists $\delta>0$ such that $|x-y|_{\mathbb{S}^1}<\delta$ 
implies that $|u(x)-u(y)|<r$ for all $x,y\in\mathbb{S}^1$, and hence that $M_r(u)\leq\frac{2\pi}{\delta}$. 

\qed

We now observe that the $r$-multiplicity, unlike length, is uniformly bounded for sequences of curves which converge uniformly.

\begin{prop} \label{convergence}  Let $u:\mathbb{S}^1\to\reals^2$ be a continuous map.  For each $r>0$ there is a constant $C=C(r,u)$
such that if $\gamma_n$ is a sequence of curves that converge uniformly to $u$, then 
$$
M_r(\gamma_n)\leq C
$$
for all $n$ such that $d(\gamma_n,u)<\frac{r}{4}$.
\end{prop}

\proof Given $r>0$ choose $\delta=\delta(r,u)>0$ such that  $|x-y|<\delta$ implies $|u(x)-u(y)|<r/2$.  Now suppose that $d(\gamma_n,u)<r/4$.  Then $|x-y|<\delta$ implies $|\gamma_n(x)-\gamma_n(y)|<r$, and using the argument from Proposition~\ref{multiplicity} it follows that $M_r(\gamma_n)\leq \frac{2\pi}{\delta}$.\qed


\subsection{A compactness property}~\label{compactness_section} We now investigate a family of curves whose $r$-multiplicity is uniformly bounded at each $r>0$.  The results of this section are not necessary for the sequel, but are included to demonstrate the nature of $r$-multiplicity. 

Let $\widetilde{\mathcal{C}}(\mathbb{S}^1,\reals^2)$ be the space of continuous maps which are locally nonconstant.  Given a function $f:\reals^+\to\reals^+$ define 
$$\mathcal{U}_f=\{u\:|\:M_r(u)<f(r) \:\:\forall \:r>0\}\subset\widetilde{\mathcal{C}}(\mathbb{S}^1,\reals^2),
$$
and let $\mathcal{U}_{f,R}$ be the subset of $\mathcal{U}_f$ whose curves are contained in $B_R(0)$.

The family $\mathcal{U}_{f,R}$ is not equicontinuous since the $r$-multiplicity does not depend on the parametrization.  However, Theorem~\ref{equicontinuity} shows that one can apply the Arzel$\grave{a}$-Ascoli Theorem, and extract convergent subsequences, once a suitable parametrization has been chosen for each map in $\mathcal{U}_{f,R}$.

\begin{thm} ~\label{equicontinuity}  For each $f:\reals^+\to\reals^+$ and $R>0$ there is a function $\delta:\reals^+~\to~\reals^+$ with the following property:  If $u\in\mathcal{U}_{f,R}$ then there exists a homeomorphism $\phi:\mathbb{S}^1\to\mathbb{S}^1$ such that $\delta$ is a modulus of continuity for $u\circ\phi$.
\end{thm}

\proof Fix a decreasing sequence $\{r_n\}$ with $r_n\to 0$ as $n\to\infty$. 

The argument consists of an iterative procedure performed at successively smaller scales.  Let $\{x_1,x_2,\ldots, x_{n_1}\}\subset\mathbb{S}^1$ be a maximal ordered set such that $d(u(x_i), u(x_{i+1}))\geq\frac{r_1}{8}$ for each $i\in\{1,2,\ldots, n_1\}$.  The maximality of $\{x_1,x_2,\ldots, x_{n_1}\}$ implies that $\diam(u[x_i,x_{i+1}])<\frac{r_1}{2}$ for each $i$.

The first step in defining the homeomorphism $\phi:\mathbb{S}^1\to\mathbb{S}^1$ is the following:  For each $i\in\{1,2,\ldots,n_1\}$ set
$$
\phi\left(\frac{2\pi i}{n_1}\right)=x_i.
$$

{\bf Claim 1:} There exists a constant $C_1=C_1(f,R,r_1)$ such that $n_1<C_1$.

\proof [Proof of Claim:] Fix the constant $c=(16\sqrt{2})^{-1}$.  Let $\Lambda=\{\ell_j\}$ be the set of horizontal and vertical lines through the grid $cr_1\integers^2$ that intersect $B_R(0)$.  Then $d(x_i,x_{i+1})\geq\frac{r_1}{8}$ implies that there exists parallel lines $\ell_{j_1},\ell_{j_2}\in\Lambda$ such that $u[x_i,x_{i+1}]\cap\ell_{j_i}\neq\emptyset$ for $i=1,2$.   Thus
$$
n_1\leq M_{cr_1}(u)|\Lambda|<f(cr_1)\left(\frac{4R}{cr_1}+2\right),
$$
and this proves the claim.

Let $\psi:\mathbb{S}^1\to\mathbb{S}^1$ be any homeomorphism extending $\phi$ and let $v=u\circ\overline{\phi}$.  Then $d_{\mathbb{S}^1}(p,q)<2\pi/C_1$ implies that $v([p,q])\subset u[x_i,x_{i+2}]$ for some $i$ and so $d(v(p),v(q))<r_1$ since $\diam(u[x_i,x_{i+1}])<\frac{r_1}{2}$ for each $i$ .  This shows that the modulus of continuity of $v$ is controlled at all scales $r\geq r_1$ by setting $\delta(r)=2\pi/C_1$ for all $r\geq r_1$. 

We then repeat the procedure at the scale $r_2$:  

For each $i\in\{1,2,\ldots,n_1\}$, let $\{x_i=y_0,y_1,\ldots, y_{n_{2i}}=x_{i+1}\}$ be a maximal ordered set of points in $[x_i,x_{i+1}]\subset\mathbb{S}^1$ such that $d(u(y_i),u(y_{i+1}))\geq\frac{r_2}{8}$.  For each $j\in\{0,1,\ldots,n_{2i}\}$ set

$$
\phi\left(\frac{2\pi i}{n_1}+\frac{2\pi j}{n_1n_{2i}}\right)=y_j.
$$

As above there exists $C_2=C_2(f,R,r_1,r_2)$ such that $n_{2i}<C_2$ for all $i$, and if $\psi$ is an extension of $\phi$ to a homeomorphism, then the modulus of continuity of $u\circ\psi$ is controlled in terms of $f$ and $R$ at scales $r\geq r_2$. 

At the end of this iterative process we obtain a map $\phi:E\to\mathbb{S}^1$, where $\overline{E}=\mathbb{S}^1$ since $u$ is not constant on an interval.  The fact that $\phi$ extends to a homeomorphism of $\mathbb{S}^1$ follows from the following two observations:  First, $u\circ\phi$ is uniformly continuous on $E$ (with modulus of continuity $\delta$).  Second, given $\epsilon>0$ there is a $\rho>0$ such that $\diam(u[x,y])<\rho$ implies $|x-y|<\epsilon$.  Indeed, if not, then there is an $\epsilon>0$ and sequences $\{y_n\}$ and $\{z_n\}$ in $\mathbb{S}^1$ so that $d(y_n,z_n)\geq\epsilon$, but $\diam(u[y_n,z_n])\to 0$ as $n\to\infty$.  After, taking convergent subsequences we get points $y_0$ and $z_0$ such that $d(y_0,z_0)\geq\epsilon$, but $\diam(u[y_0,z_0])=0$, contradicting the fact that $u$ is not constant on an interval.  

Combining these two statements we see that $\phi$ itself is uniformly continuous on $E$ and hence that it can be extended to a continuous map $\phi:\mathbb{S}^1\to\mathbb{S}^1$, which is a homeomorphism since $\phi|_E$ is strictly order-preserving.\qed


\section{The straightening lemma}\label{straight_section}

Given constants $a,b\in\reals$ and $c>0$, the 1-parameter family of
curves
$$
u(s,t)=\left(\frac{\log(\sec(cs))}{c}+ct -a, s-b\right),
\:\:\:\:\:\frac{-\pi}{2}\leq s\leq\frac{\pi}{2c}
$$
is a translating solution to CSF, which moves with speed $c$ in the
direction of the positive $x$-axis.  In this section we use these so-called grim repears to
prove Lemma~\ref{straightening_intro}, which controls the time it takes non-compact curves that are linear
outside a compact set to unwind.  We start by constructing a grim reaper where the slope of the tangent
at any point which passes through a neighbourhood of the origin is
controlled.

\begin{lem} \label{reaper} Given $r,\alpha>0$, there exist constants
$a,b$ and $c$ such that with $u(s,t)$ defined as above we have
\newline\indent(1) $(-r,r)^2$ is contained in the interior of
the convex hull of $u(s,0)$, and
\newline\indent(2) if $s-b\geq -r$, then the slope of
the tangent to the curve at $u(s,0)$ is less than $\alpha$.
\end{lem}

\proof Given $s_0\in(0,\frac{\pi}{2c})$, the slope of the tangent to
the curve $u(s,0)$ at $s=s_0$ is $\tan(cs)^{-1}$.  Thus $s\geq
c^{-1}\arctan(\alpha^{-1})$ implies that the slope of the tangent is not
more than $\alpha$.  Defining
$$
c=(3r)^{-1}\left(\frac{\pi}{2}-\arctan(\alpha^{-1})\right)
$$
and
$$
b=\frac{\pi}{2c}-2r
$$
ensures that $s\geq c^{-1}\arctan(\alpha^{-1})$ whenever $s-b\geq -r$,
and hence that (2) holds.  In addition, $c$ and $b$ have been chosen
so that $u(s,0)$ is asymptotic to the line $y=\frac{\pi}{2c}-b=2r$
as $s\to\frac{\pi}{2c}$.  This allows us to define

$$
a=r+\frac{\log(\sec(c(b+r)))}{c}
$$
so that the point $(-r,r)$ lies on $u(s,0)$.\qed

\begin{rem} \label{rem_quant}Note that $c=f(\alpha)r^{-1}$ for some function $f$
satisfying $f(\alpha)\to 0$ as $\alpha\to 0$.  Since the speed of $u(s,t)$ is
inversely proportional to the width of its opening this choice for
$c$ is, up to a constant, as large as possible.  Also $a=g(\alpha)r$, and
$g(\alpha)\to\infty$ as $\alpha\to 0$.
\end{rem}

\begin{figure}
\centering
\scalebox{0.9}{\includegraphics{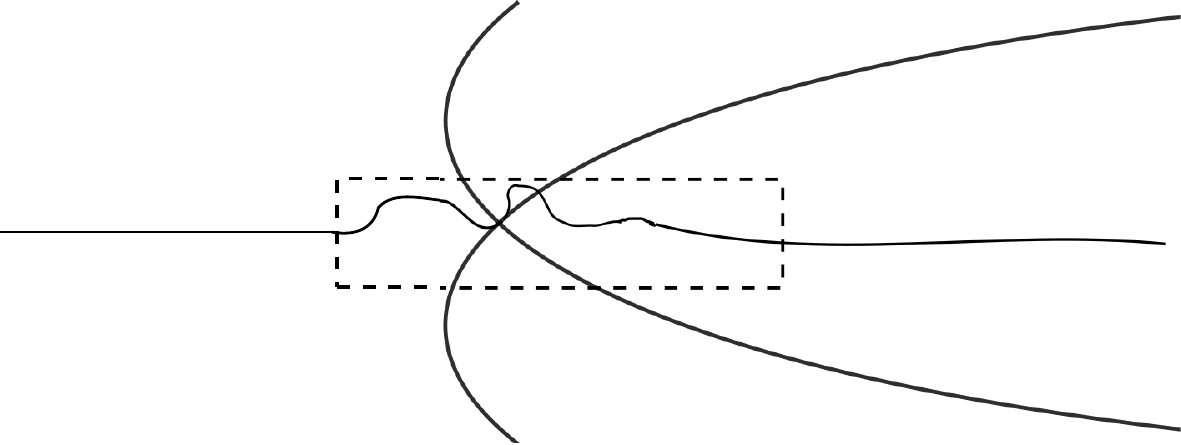}}
\caption{At time $T$, any point in $\gamma_T$ lies on the time $T$ evolution of the grim reapers $u_\lambda(s,t)$ and $v_\psi(s,t)$.  If the tangent to $\gamma_T$ at the intersection point is not bounded by slope of the tangents of these grim reapers (as in the Figure), then we obtain a contradiction since $\gamma_T$ must then intersect one of $u_\lambda(s,T)$ and $v_\psi(s,T)$ a second time.}
\label{Reaper3}
\end{figure}

\proof [Proof of Lemma~\ref{straightening_intro}]  Given $r,\alpha>0$, let $u(s,t)$ be the solution to CSF given by
Lemma~\ref{reaper}.  Condition (1) in Lemma~\ref{reaper} implies that
$\gamma$ and $u(s,0)$ intersect at a single point, and hence that
$\gamma_t$ and $u(s,t)$ intersect in a single point for all
$t\geq0$.  Note that $\gamma_t$ exists for all $t>0$ by the results of .

We will need to use not only the solution $u(s,t)$, but also
particular translates and reflections.  For
$\lambda\leq 0$, define $u_\lambda(s,t)=u(s,t)+(\lambda,0)$, and let
$v_\lambda(s,t)$ be the reflection of $u_\lambda(s,t)$ in the
$x$-axis.  Like $u(s,t)$, each of these translating curves intersect
$\gamma_t$ exactly once for all $t\geq 0$.

Define
$$
T=\frac{l+a}{c}
$$
so that at time $T$ the solution $u_0(s,t)$ has passed through $[0,l]\times[-r,r]$. 
Now suppose that $\gamma_T(s_0)$ has $x$-coordinate less than $l$.  Then there exists unique $\lambda$ and $\psi$ such that
$\gamma_T(s_0)$ lies on $u_\lambda(s,T)$ and $v_\psi(s,T)$.  We
claim that the tangent to $\gamma_T$ at $\gamma_T(s_0)$ must have
slope between $-\alpha$ and $\alpha$.  Indeed, otherwise $\gamma_T$ would enter
a region contained in the convex hull of exactly one of
$u_\lambda(s,T)$ and $v_\psi(s,T)$.  See Figure~\ref{Reaper3}.  Since $\gamma_T$ eventually
lies in the intersection of the convex hulls, this implies that it
must intersect one of the two curves a second time, and this is a
contradiction.  When the $x$-coordinate is not less than $x$ one can
use grim reapers traveling in the opposite direction.

Remark~\ref{rem_quant} shows that $T(r,l,\alpha)=f(\alpha)^{-1}r(l+rg(\alpha))$,
which proves the final statement.\qed


\section{Separating foliations}~\label{holonomy_section}

The main result of this section is a different type of straightening lemma.  Given a foliation of an open ball by smooth $\alpha$-Lipschitz graphs, the bilipschitz 
constant needed to straighten the foliation to one with linear leaves can be arbitrarily large.  We prove that for any $t>0$, CSF separates the leaves uniformly, so that a bound on the bilipschitz constant can be obtained.  

For any family of smooth curves $\mathcal{F}$,
we define $\mathcal{F}_t=\{\gamma_t\mid\gamma\in \mathcal{F}\}$ to be the time $t$
evolution of the family under CSF, provided each $\gamma_t$ exists.

\begin{lem} ~\label{holonomy} Given $0<\alpha<1$ and $t, r>0$, there is a constant
$\widehat{C}=\widehat{C}(t,\alpha,r)$ such that the following holds:  Let $\ell$ and $\ell'$ be perpendicular lines whose intersection point is $x$.
Let $\mathcal{F}$ (resp.~$\mathcal{F}'$) be a smooth foliation of $\mathcal{N}_r(\ell)$  (resp.~$\mathcal{N}_r(\ell')$), whose leaves are entire
$\alpha$-Lipschitz graphs over $\ell$ (resp.~$\ell'$).
Then there is an open set $V\subset\reals^2$ and a $\hat{C}$-bilipschitz diffeomorphism $\Phi:B_r(x)\to V$
which sends the leaves of $\mathcal{F}_t$ and $\mathcal{F}'_t$ to horizontal and
vertical segments.
\end{lem}

\begin{rem}  (1) Here, and in Section~\ref{length_estimates_section}, it is not necessary for the foliations to be perpendicular.  One can take any two transverse directions, as long as $\alpha$ is chosen sufficiently small.

(2) In Section~\ref{length_estimates_section}, we obtain bounds on the length of an evolving curve by controlling the number of times a curve intersects such foliations.  
To demonstrate the difference in behaviour between Lipschitz foliations and the standard grid
we note that for any $\epsilon>0$, there is a pair of foliations of $[0,1]^2$ whose leaves are $\epsilon$-close to
horizontal (resp. vertical) lines in the $C^{1/\epsilon}$-topology, and a collection of curves $\{\gamma_i\}$ such that
each leaf is hit precisely once, but $\Sigma_i \mathscr{L}(\gamma_i)> 1/\epsilon$.  
\end{rem}

To prove Lemma~\ref{holonomy} it suffices to show that the
holonomy maps of $\mathcal{F}_t$ and $\mathcal{F}'_t$ have bounded derivative.  Recall
that given a foliation $F$, and a set of transversals $\{T_i\}$ the
holonomy maps $\Phi$ are defined as follows:  Given $x_1$ and $x_2$
in the same leaf of $F$, let $T_1$ and $T_2$ be the transversals
through $x_1$ and $x_2$ respectively.  For each $x\in T_1$, let
$\ell_x$ be the leaf of $F$ containing $x$.  Then $\Phi:T_1\to T_2$
is the locally defined diffeomorphism given by $\Phi(x)=\ell_x\cap
T_2$.  In this example we may take the leaves of one foliations as transversals for the other.

A main ingredient in the proof of Lemma~\ref{holonomy} is the linearization of the curve shortening equation, which controls the derivative of the holonomy maps.  We obtain the necessary bounds using the Harnack Inequality, which we state now for convenience.

\begin{thm} [Harnack Inequality for linear parabolic PDE's ~\cite{EV}]~\label{harnack} Suppose $Lu=au_{xx}+bu_x+cu$ and that $u(x,t)\in C^2(U\times (0,T))$ solves
$$
u_t-Lu=0,
$$
and that $u(x,t)\geq 0$ on $U\times(0,T)$.  Then for any compact set $K\subset U$, and $0<t_1<t_2\leq T$, there exists a constant $C$ such that
$$
\sup_Ku(\cdot,t_1)\leq C\inf_Ku(\cdot,t_2).
$$

The constant $C$ depends on $K, U, t_1, t_2$ and the coefficients of $L$.
\end{thm}

We will apply the Harnack Inequality to the solution of an operator $Lu=au_{xx}+bu_x+cu$ with $c=0$.  In this case, sup$_x\{u(x,t)\}$ is decreasing and so the conclusion 
holds with $t_1=t_2$.

\proof [Proof of Lemma~\ref{holonomy}] We may assume that $\ell$ is parallel to the $x$-axis.  We first note that since each curve in $\mathcal{F}$ has a unique solution to CSF that exists for all time, that $\mathcal{F}_t$ continues to foliate the strip $\mathcal{N}_r(\ell)$.  Given $t>0$, let $\Psi:T_1\to T_2$ be a holonomy map of the foliation $\mathcal{F}_t$ restricted to $B_r(x)$, and choose $x_1\in T_1$.    Let  $u(x,t)$ be the leaf of $\mathcal{F}_t$ containing $x_1$, written as a graph over the $x$-axis.  Fix $0<\mu<1$.   According to~\cite{EH91}, the curvature of $u(x,t)$ is bounded by a constant $C'=C'(\alpha,t)$ on $[\mu t,t]$.  

The graphical form of the curve shortening flow equation is $u_t=\frac{u_{xx}}{1+u_x^2}$, and the linearization of CSF at $u(x,t)$ is

$$
w_t=\frac{w_{xx}}{1+u_x^2}-\frac{2u_xu_{xx}}{(1+u_x^2)^2}w_x.
$$

The coefficient of $w_{xx}$ is well controlled since $|u_x|<\alpha$, and for the coefficient of $w_x$ on $[\mu t, t]$ we have

$$
\left\arrowvert\frac{2u_xu_{xx}}{(1+u_x^2)^2}\right|<\frac{2\alpha |u_{xx}|}{(1+u_x^2)^{\frac{3}{2}}}=2\alpha|\kappa(x,t)|<2\alpha C'.
$$

Let $u_\delta$ be a parametrization of the leaves of $\mathcal{F}_t$ with $x_1\in u_0$.  Each $u_\delta(x,t)$ is the evolution of a particular leaf written as a graph over the x-axis.  Define

$$
v(x,\mu t)=\left |\frac{du_\delta}{d\delta}(x,\mu t)|_{\delta=0}\:^\perp\right |,
$$
where $\perp$ indicates the projection onto the normal of $u_0(x,\mu t)$, and note that $v(x,\mu t)$ exists since $\mathcal{F}_0$ is smooth and $\mathcal{F}_t$ has been evolving by CSF.  In addition, $v(x,\mu t)>0$ since $\mathcal{F}_{\mu t}$ is a foliation.   Let $v(x,t)$ be the solution of the linearized equation above with initial condition given by $v(x,\mu t)$. 

By the Harnack Inequality there is a constant $\widetilde{C}=\widetilde{C}(t,r)$ such that 

$$
\sup_K v(\cdot,t)<\widetilde{C}\inf_K v(\cdot,t),
$$
where $K$ is an appropriately chosen compact set whose size depends only on~$r$.  
This proves that $|D\Psi(x_1)|$ is bounded in terms of $\alpha$, $t$, and $r$ since there is a constant $c=c(\alpha)>1$ such that 
$$
c^{-1}\frac{v(x_2,t)}{v(x_1,t)}\leq |D\Psi(x_1)|\leq c\frac{v(x_2,t)}{v(x_1,t)}.
$$
\qed


\section{Constructing the initial foliation}\label{construct_section}

In this section we construct the initial foliations that ultimately lead to the length estimate in the next section.  To motivate Lemma~\ref{foliate1} we refer the reader to Section~\ref{intro_outline} for an outline of the proof of main length estimate.

The crucial property of the foliations constructed in Lemma~\ref{foliate1} is that the number of times each leaf intersects the given curve $\gamma$ is bounded by the $r$-multiplicity and the number of intersections of $\gamma$.  In order for this to work we need to count intersections with multiplicity:

\begin{defn}  [Intersection number] \label{intersection} Let $\Lambda=\{\gamma_i\}$ be a collection of arcs.  Define $\mathcal{I}(\Lambda)$ to be the number of points in the preimage of $\Lambda$ whose image is not unique.  In particular, for a closed curve $\gamma$, 
$$
\mathcal{I}(\gamma)=|\{x\in\mathbb{S}^1\:|\:\gamma(x)=\gamma(y)\:\mathrm{for\: some}\: y\in\mathbb{S}\setminus x\}|.
$$ 
Furthermore, if $\phi$ is an embedded curve then define $\mathcal{I}(\gamma,\phi)$ to be the number of points in the preimage of $\gamma$ which map into $\phi$. 
\end{defn}

We can now state the main result of this section.  Let $K$ be the minimal simply-connected compact set containing $\mathcal{N}_1(\gamma)$.

\begin{thm} \label{foliate1} Let $\gamma$ be a smooth immersed closed curve, $\ell$ be a line and
$r>0$.  Then there exists a 1-parameter family of smooth curves
$\mathcal{F}=\{\ell_x\}_{x\in[0,1]}$ such that:
\newline\indent (1) $\mathcal{F}$ foliates a region containing
$\mathcal{N}_r(\ell)$,
\newline\indent (2) $\ell_x\setminus K$ is contained in a line parallel to $\ell$ for each $x\in[0,1]$,  and
\newline\indent (3) $\mathcal{I}(\gamma,\ell_x)\leq 2(\widetilde{M}_{r,\ell}(\gamma)+\mathcal{I}(\gamma))$ for
each $x\in [0,1]$.
\end{thm}

The following lemma allows us to construct the foliation by taking the image of parallel lines under a diffeomorphism.  Recall that a proper arc in a surface $D$ is an arc whose intersection with $\partial D$ consists exactly of its endpoints. 

\begin{lem}\label{isotopy_lemma} Let $D\subset\reals^2$ be the closed region between two distinct parallel lines $\ell_0$ and~$\ell_1$.  Let $\Lambda=\{\gamma_i\}$ be  a finite collection of smooth proper arcs in $D$ which intersect $\partial D$ transversally.  Then there is a diffeomorphism $\Phi:D\to D$ supported in the interior of $D$ such that
$$
\sum_i \mathcal{I}(\Phi(\gamma_i),\ell)\leq 2(|\Lambda|+\mathcal{I}(\Lambda))
$$
for each line $\ell$ in $D$ parallel to $\ell_0$.
\end{lem}


In particular, there exists a foliation of $D$ by smooth complete curves that intersect $\cup\gamma_i$ at most $2(|\Lambda|+\mathcal{I})$ times. The proof is by induction on the number of intersection points in $\cup_i\gamma_i$.  Only the base case (when $\mathcal{I}=0$) is needed to prove Theorem~\ref{main_estimate} for embedded curves.

\proof  First suppose that $\mathcal{I}=0$.  We may assume that $\ell_0$ is parallel to the $x$-axis.  Since $\gamma_i$ intersects transversally, there exists a set of smooth arcs $\{\rho_i\}$ and $\epsilon>0$ such that
\newline\indent (ii) $\rho_i\cap\mathcal{N_\epsilon}(\partial D)=\gamma_i\cap\mathcal{N_\epsilon}(\partial D)$, and
\newline\indent (iii) $\rho_i$ has a single local minimum or maximum if the endpoints of $\gamma_i$ are contained in the same component of $\partial D$ and otherwise $\rho_i$ contains no local maximum or minimum.


The arc $\gamma_1$ is smoothly isotopic to $\rho_1$ and using basic facts from differential topology it
 extends to an isotopy of $D$ which is the identity in a neighbourhood of $\partial D$.  Thus we obtain a diffeomorphism $\Upsilon_1:D\to D$ such that
$\Upsilon_1(\gamma_1)=\rho_1$.  We then repeat this process inductively, noting that
$\Upsilon_n(\gamma_{n+1})$ is smoothly isotopic to
$\rho_{n+1}$ by an isotopy which is supported in
$D\setminus\cup_{i=1}^n\rho_i$.  After $|\Lambda|$
steps, we obtain the required  diffeomorphism $\Upsilon=\Upsilon_{|\Lambda|}$.

Now suppose $\mathcal{I}>0$.  Let $\ell_2$ be a line parallel to $\ell_0$ in the interior of $D$, and let $D_0$ (resp. $D_1$) be the closed region between $\ell_2$ and $\ell_0$ (resp. $\ell_1$).  By choosing an initial diffeomorphism we may assume that


\indent (i) there is a single point $x\in D_0$ where intersections of $\cup\gamma_i$ occur, 
\newline\indent (ii)  each individual arc in $(\cup\gamma_i)\cap D_0$ is embedded,
\newline\indent (iii) each arc $\gamma_i$ intersects $\ell_2$ transversally, and 
\newline\indent (iv) any arc of $(\cup\gamma_i)\cap D_0$ with both its endpoints in $\ell_2$ passes through $x$.

\smallskip

In order to complete the proof we need to count the arcs in $D_0$ according to their endpoints.  Let $\chi$ be the set of arcs in $(\cup\gamma_i)\cap D_0$, and let $A$ be the number of arcs in $\chi$ that have both endpoints in $\ell_0$, and pass through $x$.
Let $B$ be the number of arcs in $\chi$ that have one endpoint in each of $\ell_2$ and $\ell_0$, and pass through $x$.
Let $C$ be the number of arcs in $\chi$ that have both endpoints in $\ell_2$, and pass through~$x$.

Let $\alpha$ be the number of arcs in $\chi$ that have both endpoints in $\ell_0$, but do not pass through $x$.
Let $\beta$ be the number of arcs in $\chi$ that have one endpoint in each of $\ell_2$ and $\ell_0$, but do not pass through $x$.
By assumption (iv) above there are no arcs in $\chi$ disjoint from $x$ with both endpoints in $\ell_2$.

Then the number of arcs in $D_1$ is  $|\Delta|+C-A-\alpha$ and the number of intersection points in $D_1$ is $\mathcal{I}-(A+B+C)<\mathcal{I}$.  Thus by induction there exists a diffeomorphism $\Phi_1:D\to D$ supported in the interior of $D_1$ such that 
\begin{align*}
\sum_i \mathcal{I}(\Phi_1(\gamma_i),\ell)&\leq 2(|\Delta|+C-A-\alpha+\mathcal{I}-(A+B+C))
\\ &\leq 2(|\Delta|+\mathcal{I}).
\end{align*}
for each line $\ell$ between $\ell_1$ and $\ell_2$.

In $D_0$ each arc intersects $x$ at most once, and so we can use a procedure as in the base case to isotope each component to an arc with at most one local maximum or minimum depending on the endpoints of the arc.  Thus, there is a diffeomorphism $\Phi_0:D\to D$ supported in the interior of $D_0$ such that 
\begin{align*}
\sum_i \mathcal{I}(\Phi_0(\gamma_i),\ell)&\leq (2\alpha+\beta)+(2A+B+2C)
\\ &\leq 2(|\Delta|+\mathcal{I}).
\end{align*}
for each line $\ell$ between $\ell_0$ and $\ell_2$.  Since $\Phi_0$ and $\Phi_1$ are supported on disjoint sets the required diffeomorphism is obtained by composing.\qed

\begin{figure}
\centering
\scalebox{0.8}{\includegraphics{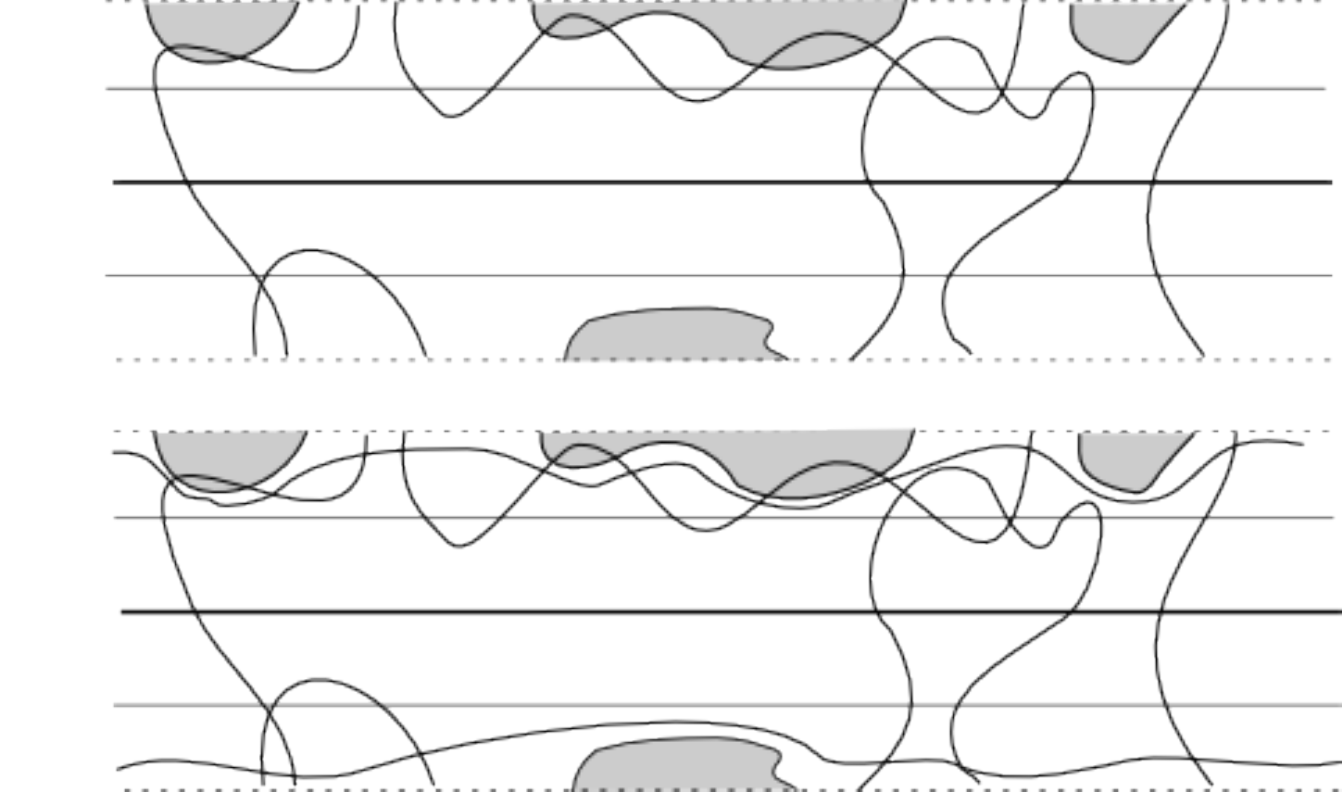}}
\caption{In the top picture is an example where $\widetilde{M}=5$.  Two of the components defining $\widetilde{M}$ split into separate arcs in $\Omega$, the unshaded region.  The second picture includes the outermost leaves of $\mathcal{F}$.  The leaf running along the upper boundary of $\Omega$ intersects $\gamma$ 10 times, while the leaf running along the lower boundary of $\Omega$ intersects $\gamma$ 6 times.}  
\label{imm_fol_fig}
\end{figure}

We are now able to construct the initial foliation.  

\proof  [Proof of Theorem~\ref{foliate1}] Let $\{\gamma_i\}_{i\in\Delta}$ be the arcs in $\gamma$ defining $\widetilde{M}_{r,\ell}(\gamma)$.
See Definition~\ref{intro_mult_def}.

The first step is defining the outermost leaves of $\mathcal{F}$.  Let $\varphi=\gamma\setminus\cup\gamma_i$.  Define $\Omega\subset\mathcal{N}_{2r}(\ell)$ to be the connected component of $\mathcal{N}_{2r}(\ell)\setminus\varphi$ containing $\ell$, and let $\Omega_1$ be one of the two components of $\partial\Omega$.  For $\epsilon>0$ sufficiently small the intersection of $\gamma$ with the open strip $S=\Omega\cap\mathcal{N}_\epsilon(\Omega_1)\subset\mathcal{N}_{2r}(\ell)\setminus\mathcal{N}_r(\ell)$ consists of embedded arcs which cross the strip since there are only a finite number of intersections.  Thus there exists $\ell_0\in S$ satisfying:

\indent(i) $\ell_0\setminus K$ is contained in a line,
\newline\indent (ii) $\ell_0$ intersects $\gamma$ transversally, 
\newline\indent (iii) $\ell_0$ intersects $\gamma$ only in $\cup_i\gamma_i$, and 
\newline\indent (iv) for each $i\in\Delta$, $\ell_0$ intersects each component $\rho$ of $\gamma_i\cap\Omega$ exactly once for each for endpoint of $\rho$ that lies in $\Omega_1$.

See Figure~\ref{imm_fol_fig}.

We define $\ell_1$ similarly, except lying near the other component of $\partial\Omega$.   

For each $i\in\Delta$, let $\Gamma_i$ be the number of components of $\gamma_i\cap\Omega$, and let $D$ be the closed region between $\ell_0$ and $\ell_1$.  Then $D\cap\gamma$ consists of $\sum \Gamma_i$ proper arcs, and since $\gamma_i$ intersects $\varphi$ at least $\Gamma_i-1$ times, the total number of intersections in $D$ is at most $\mathcal{I}(\gamma)-(\sum\Gamma_i-1)$.  By Lemma~\ref{isotopy_lemma} there exists a foliation of $D$ by smooth complete curves such that $\gamma$ hits each leaf at most 
$$
2\left(\sum_i \Gamma_i+\mathcal{I}(\gamma)-\sum_i(\Gamma_i-1)\right)=2\left(|\Delta|+\mathcal{I}(\gamma)\right)
$$
times, and it is clear that this foliation may be chosen so that the leaves are linear outside $K$.

\qed


\section{The length estimate}\label{length_estimates_section}

In this section we combine the results of the previous sections to show that the $r$-multiplicity can be used to bound the length of an evolving curve.
While Theorem~\ref{main_estimate_intro} was stated only for embedded curves the statement below includes immersed curves, with bounds depending only on the $r$-multiplicity and number of self-intersections.  Recall that $\mathcal{I}(\gamma)$ has been defined in Section~\ref{construct_section}.

\begin{thm} \label{main_estimate}
For all $t,d>0$ there exists constants $r(t,d)$ and $C=C(t,d)$ with the following property:  Let
$\gamma$ be a smooth immersed closed curve with $\diam(\gamma)<~d$.  Then for all $t>0$ such that $\gamma_t$ exists
\newline\indent(1) 
$$
\mathscr{L}(\gamma_t)<C\left(M_r(\gamma) +\mathcal{I}(\gamma)\right),
$$
and
\newline\indent(2) for any $0<\nu\leq 1$ and $x\in\reals^2$
$$
\mathscr{L}(\gamma_t\cap B_{\nu r}(x))<\nu C\left(M_r(\gamma)+\mathcal{I}(\gamma)\right).
$$
\end{thm}

Note that Corollary~\ref{convergence_estimate_intro} is an immediate consequence of Theorem~\ref{main_estimate} (1) and Proposition~\ref{convergence}.

From this point on, we fix a constant
$0<\alpha\ll 1$.  All subsequent constants depend on this choice, but
there is no need to vary $\alpha$, so this dependence is suppressed.  We begin by making precise the effect of Lemma~\ref{straightening_intro} on the foliation constructed in Section~\ref{construct_section}.

\begin{lem} \label{foliate2} Given a line $\ell$, $r>0$ and a smooth curve $\gamma$,
let $\mathcal{F}$ be a family of curves constructed in Lemma~\ref{foliate1}.
Define $t_0=T(2r,d+2,\alpha)$ where $\diam(\gamma)~<~d$ and $T$ is the function defined in Lemma~\ref{straightening_intro}.  Then for $t\geq
t_0$ such that $\gamma_t$ exists we have:
\newline\indent (1) $\mathcal{F}_t$ foliates a region which contains
$\mathcal{N}_r(\ell)$,
\newline\indent (2) each curve in $\mathcal{F}_t$ is an $\alpha$-Lipschitz graph
over $\ell$, and
\newline\indent (3) $\mathcal{I}(\gamma_t,(\ell_x)_t)\leq 2(\widetilde{M}_{r,\ell}(\gamma)+\mathcal{I}(\gamma))$ for
each $x\in [0,1]$.
\end{lem}

\proof The region between $({\ell_0})_t$ and $({\ell_1})_t$ is foliated
for each $t>0$, and (1) follows since $\partial \mathcal{N}_r(\ell)$ acts as a
barrier for the evolutions of $\ell_0$ and $\ell_1$.

The second condition follows directly from Lemma~\ref{straightening_intro},
and the choice of $t_0$, while~(3) holds since the number of
intersections does not increase along the flow~\cite{An88}, even for immersed curves.\qed

We now complete the proof of Theorem~\ref{main_estimate}.

\proof [Proof of Theorem~\ref{main_estimate}] Let $\gamma$ be a
smooth curve with $\diam(\gamma)<d$.  Given $t>0$, define
the scale $r$ so that $t=2T(2r,d+2,\alpha)$, where $T$ is the function in
Lemma~\ref{straightening_intro}.  The extra factor of 2 is included since we plan to straighten the leaves of $\mathcal{F}$ to 
$\alpha$-Lipschitz graphs during the interval $[0,t/2]$, and then apply Lemma~\ref{holonomy} during $[t/2,t]$. 

Given $x\in\reals^2$, let $\ell$ and $\ell'$ be two perpendicular lines through $x$.  Applying 
Lemma~\ref{foliate1} with $\ell$ (resp. $\ell'$) we obtain a foliation $\mathcal{F}$ (resp.
$\mathcal{F}'$) of $\mathcal{N}_r(\ell)$ (resp. $\mathcal{N}_r(\ell')$). 
Let $0<\nu\leq 1$.  By Lemma~\ref{foliate2}, $\mathcal{F}_{t/2}$ and $\mathcal{F}'_{t/2}$ foliate $B=B_{\nu r}(x)$, and each leaf is an
$\alpha$-Lipschitz graph over the appropriate axis.  Let
$\widehat{C}=\widehat{C}(t/2,\alpha,r)$ be the constant defined by
Lemma~\ref{holonomy}.  Then there exists a $\widehat{C}$-bilipschitz
diffeomorphism~$\Phi$, that straightens the leaves of $\mathcal{F}_t$ and
$\mathcal{F}'_t$.  Then $\Phi(\gamma_t\cap B)$ is contained in a ball of
radius $\widehat{C}\nu r$ and intersects each vertical and horizontal
line at most $2(\widetilde{M}_r(\gamma)+\mathcal{I}(\gamma))\leq 2(M_r(\gamma)+\mathcal{I}(\gamma))$ times.  Then (2) follows from
$$
\mathscr{L}(\Phi(\gamma_t\cap B))< 8\nu r\widehat{C}(M_r(\gamma)+\mathcal{I}(\gamma))
$$
which implies
$$
\mathscr{L}(\gamma_t\cap B)< 8\nu r\widehat{C}^2(M_r(\gamma)+\mathcal{I}(\gamma)).
$$

But then at most $\left(\frac{d}{r}\right)^2$ balls of radius $r$ are needed to
cover $\gamma_t$, and so
$$
\mathscr{L}(\gamma_t)<\frac{8d^2\widehat{C}^2}{r}(M_r(\gamma)+\mathcal{I}(\gamma)),
$$
and (1) holds.

\qed

Besides the length bounds for an approximating sequence given by Corollary~\ref{convergence_estimate_intro} there is also the following curvature bound:

\begin{cor} \label{curvature_estimate} Let $u:\mathbb{S}^1\to\reals^2$ be a continuous map, and let $\gamma_n$ be a sequence of embedded smooth closed curves converging uniformly to $u$.  Then for any $\epsilon>0$ there exists $0<t<\epsilon$ such that
$$
\liminf_{n\to\infty} \int_{{\gamma_n}_t}\kappa^2<\infty.
$$

\end{cor}

\proof  We may assume that $\epsilon>0$ is small enough so that $({\gamma_n})_\epsilon\neq\emptyset$ for sufficiently large $n$. 

By Theorem \ref{main_estimate} there exists a constant $C$ such that  $\mathscr{L}({\gamma_n}_{\frac{\epsilon}{2}})<C$ for large~$n$.  Consider the functions
$$
f_n(u)=\int_{\gamma_{n_u}}\kappa^2ds
$$
defined for $u\in[\epsilon/2,\epsilon]$.  Since $-f_n(u)$ is the derivative of $\mathscr{L}({\gamma_n}_u)$ under
CSF, we obtain
$$
\int_{\epsilon/2}^\epsilon f_n(u)du<C.
$$
But then $\{f_n\}$ is a collection of continuous functions on a
compact interval, with bounded $L^1$-norm, and so the set $\{u\mid
$ liminf$_nf_n(u)=\infty\}\subset [\epsilon/2,\epsilon]$ has Lebesgue measure zero.\qed  


\section{Locally-connected sets} ~\label{lc_sets_section}

In order to prove Theorem~\ref{main} we approximate a compact set $K$ by smooth curves in $\reals^2\setminus K$.  We refer the reader to Section~\ref{LSF_section} for a description of the level set flow in terms of such an approximation.  In this section we prove that in each component of the complement, a sequence of approximating curves can be chosen 
so that the $r$-multiplicity is uniformly bounded for each $r>0$.  The bounds on $r$-multiplicity depend on the local-connectivity function of $K$, which we define now.

Let $K\subset\reals^2$ be locally-connected and compact.  Given $p,q\in K$, let $\diam(p,q)$ be the minimum diameter of a connected subset of $K$ containing both $p$ and $q$, and define the local-connectivity function $f:\reals^+\to\reals^+$ of $K$ by 

$$
f(s)=\inf\{r\mid d(p,q)\leq s\Rightarrow \diam(p,q)\leq r\}.
$$
\begin{prop}\label{local_function} Let $K\subset\reals^2$ be locally-connected and compact, and let $f$ be the local-connectivity function of $K$.  Then 
\newline\indent (1) $f$ is non-decreasing,
\newline\indent (2) $f(s)\geq s$ for each $0<s\leq\diam(K)$, and 
\newline\indent (3) $\lim_{s\to 0^+}f(s)=0$.
\end{prop}

\proof  The first two properties follow straight from the definition.  For (3), suppose that there exists $r>0$ such that $f(s)>r$ for all $s>0$.  Then there exists sequences $p_n$ and $q_n$ such that $d(p_n,q_n)\to 0$, but $\diam(p_n,q_n)>r$.  After taking a convergent subsequence we may assume $p_n,q_n\to x_0$ for some $x_0\in K$, and since the assumptions on $K$ imply that it is locally-path-connected at $x_0$ this leads to a contradiction.\qed

\begin{thm} \label{approx_bound} Let $K\subset\reals^2$ be locally-connected, connected and compact, and let $\Gamma$ be a component of $\reals^2\setminus K$.  Then there exists a sequence $\{D_n\}$ of sets such that 
\newline\indent (1) $\partial D_n$ is a closed curve,
\newline\indent (2) $D_n\subset D_{n+1}$ for each $n$,
\newline\indent (3) $\Gamma=\cup_n D_n$, and
\newline\indent (4) for each $r>0$, there exists $\widetilde{C}=\widetilde{C}(r,K)$ such that $M_r(\partial D_n)<\widetilde{C}$.
\end{thm}

\begin{rem} ~\label{rem_approx} In the proof below each curve $\partial D_n$ is a polygon.  In order to approximate by smooth curves and obtain the same conclusions as in Theorem~\ref{approx_bound}, one can choose a smooth curve sufficiently close to $\partial D_n$  making sure not to increase the $r$-multiplicity.  This can be done, for example, by evolving the polygon by CSF for a short time.  We note that if $\gamma_n$ is such a smooth sequence, then the conclusions of both Corollary~\ref{convergence_estimate_intro} and Corollary~\ref{curvature_estimate} hold. 
\end{rem}

\proof  Fix $x\in\Gamma$.  For each $i\in\naturals$, consider the square grid in $\reals^2$ whose vertices are $\{(\frac{i}{2^n},\frac{j}{2^n})\mid i,j\in\integers\}$.   Define $D_n\subset\Gamma$ by including all closed squares of this grid that can be joined to the square containing $x$ by a sequence of closed squares such that \newline\indent (1) each closed square lies entirely within $\Gamma$, and \newline\indent (2) successive squares share a common side. 

If $D_n\neq\emptyset$, then $\partial D_n$ is an embedded polygon since $K$ is connected.  Also, the sequence $\{D_n\}$ defined in this way is clearly a nested exhaustion of $\Gamma$.  

Fix $r>0$.  By Proposition~\ref{local_function}(3) there exists $s>0$ such that $f(s)<\frac{r}{2}$, where $f$ is the local-connectivity function for $K$.  

Given a line $\ell$, let $\{U_i\}_{i\in\Lambda}$ be the components of $\mathcal{N}_{r/2}(\ell)\setminus D_n$ such that $\partial U_i$ contains an arc between distinct components of $\partial\mathcal{N}_{r/2}(\ell)$.  A simple counting argument gives
$$
M_{r,\ell}(\partial D_n)\leq 2|\Lambda|,
$$
with equality when $D_n$ is an unbounded region.

The construction of $D_n$ implies that if $x\in\partial D_n$, then there exists $p\in K$ such that $d(x,p)\leq\frac{\sqrt{2}}{2^n}$, and the line segment joining $x$ and $p$ intersects $D_n$ only at $x$.  Choose $n_0\in\naturals$ such that $\frac{\sqrt{2}}{2^{n_0}}<\frac{r}{2}$.

\smallskip

{\bf Claim:} If $n\geq n_0$ then $K\cap U_i\cap\ell\neq\emptyset$.

Let $\rho$ be a component of $U_i\cap\ell$ such that both paths in $\partial D_n$ between endpoints of $\rho$ exit $\mathcal{N}_{r/2}(\ell)$.  If $\rho\cap K=\emptyset$, then $K\cap U_i$ must lie entirely in one component of  $U_i\setminus\rho$ since otherwise $D_n\cup\rho$ separates $K$.  

Let $\sigma$ be the path in $\partial D_n$ between endpoints of $\rho$ such that the region bounded by $\rho\cup\sigma$ does not contain $K$.  Since $\sigma$ exits $\mathcal{N}_{r/2}(\ell)$ there exists $x\in\sigma$ such that $d(x,\rho)\geq\frac{r}{2}$.  But then any line segment beginning at $x$ which does not enter $D_n$ must travel at least $\frac{r}{2}>\frac{\sqrt{2}}{2^n}$ before intersecting $K$, and this is a contradiction.  This proves the claim.

Let $\{p_i\}$ be a collection of points such that $p_i\in K\cap U_i\cap\ell$.  Then $i\neq j$ implies that $d(p_i,p_j)>s$ since otherwise the definition of $s$ implies that there is a connected subset of $K$ containing both $p_i$ and $p_j$ with diameter less than~$\frac{r}{2}$, a contradiction.

Thus 
$$
|\Lambda|\leq2\left(\frac{\diam(K)}{s}+1\right),
$$
and the result follows since the right-hand side depends only on $K$ and $r$.\qed


\section{Level set flow in $\reals^2$}\label{LSF_section}

In this section we recall a definition of the level set flow, and give an alternate
characterization of it in the plane.  This case is
much simpler than in higher dimensions since there is only one type
of singularity for embedded closed curves evolving by CSF.  We refer the reader to~\cite{CCG91} and~\cite{ES91}  for the analytic origins of level set methods in mean curvature flow, and to~\cite{I93},~\cite{I94}, ~\cite{W00} and ~\cite{W03} for geometric treatments which developed the definition below.

\begin{defn} [Weak set flow] Let $K\subset\reals^{n+1}$ be compact,
and let $\{K_t\}_{t\geq 0}$ be a 1-parameter family of compact sets
with $K_0=K$, such that the space-time track
$\cup(K_t\times\{t\})\subset\reals^2$ is closed. Then
$\{K_t\}_{t\geq 0}$ is a {\it weak set flow} for $K$ if for every
smooth mean curvature flow $\Sigma_t$ defined on
$[a,b]\subset[0,\infty]$ we have
$$
K_a\cap\Sigma_a=\emptyset\Longrightarrow K_t\cap\Sigma_t=\emptyset
$$
for each $t\in[a,b]$.
\end{defn}

Among all weak set flows there is one which is distinguished:

\begin{defn}[Level set flow] The {\it level set flow} of a compact
set $K\subset\reals^{n+1}$ is the maximal weak set flow.  That is, a
weak set flow $K_t$ such that if $\widehat{K}_t$ is any other weak
set flow, then $\widehat{K}_t\subset K_t$ for all $t\geq 0$.
\end{defn}

The existence of the level set flow is verified by taking the
closure of the union of all weak set flows. The rest of this section is devoted to giving
an explicit description of the level set flow in $\reals^2$.

Let $K\subset\reals^2$ be compact and connected.  Then
$\reals^2\setminus K$ consists of one unbounded
component~$\Gamma_0$, and perhaps (infinitely many) other components
$\{\Gamma_i\}_{i\geq 1}$, that are open, bounded and
simply-connected.  Using CSF, we can evolve such a bounded, 
simply-connected domain $\Gamma$ in a natural way: Let $D_n$ be a
nested exhaustion of $\Gamma$ by smooth closed 2-disks.  Let
$\partial {D_n}_t$ be the evolution of $\partial D_n$ by CSF, and
define ${D_n}_t$ to be the region bounded by $\partial {D_n}_t$.
Note that ${D_n}_t$ defined in this way is simply the level set flow
of $D_n$.  Now define $\Gamma_t:=\cup_n {D_n}_t$, and note that
$\Gamma_t$ does not depend on the original exhaustion.  

If $\Gamma$ is the unbounded component of $\reals^2\setminus K$, one can define $\Gamma_t$ in a similar manner.  
In this case, there is the possibility that $m(K)=0$ 
and $\reals^2\setminus K=\Gamma$, when $\Gamma_t=\reals^2$ for all $t>0$.   

We claim that the evolutions of $\Gamma_i$, which
are independent of each other, completely determine the level set
flow of $K$.  Indeed, it is easy to check the following:

\begin{prop} ~\label{LSF} $\reals^2\setminus \bigcup_{i\geq 0}(\Gamma_i)_t$ is the level set flow of
$K$.
\end{prop}

In particular, if $\reals^2\setminus K$ contains $N$ components,
then the level set flow of $K$ is completely determined by $N$
sequences of smooth curves approaching $K$.

Since the area contained in a smooth curve decreases at a constant rate the description of the level set flow in
Proposition~\ref{LSF} allows one to compute the rate of change of $m(K_t)$ at any positive time, where $m$ denotes the Lebesgue measure.  However, at a time when a component of
the complement vanishes, the derivative from the left and the right
will not coincide.  More precisely:

\begin{prop} \label{area} Let $K\subset\reals^2$ be compact and connected, and
let $\{\Gamma_i\}_{i\geq 0}$ be the components of $\reals^2\setminus
K$.  Given $T>0$ so that $K_T\neq\emptyset$ define 
$N_T=|\{i\mid m(\Gamma_i)\}|\geq 2\pi T$ and $M_T=|\{i\mid m(\Gamma_i)\}|> 2\pi T$.  Then $N_T,M_T<\infty$ and
$$
\frac{d^-m(K_t)}{dt}|_{t=T}=2\pi(N_T-2),
$$
and
$$
\frac{d^+m(K_t)}{dt}|_{t=T}=2\pi(M_T-2),
$$
where $d^-$ and $d^+$ are the derivatives from the left and the
right respectively.
\end{prop}

If $m(K)>0$ or $\reals^2\setminus K$ consists of more than two components, then Proposition~\ref{area} implies that $m(K_t)>0$ for small $t$.
In addition, if $m(K)=0$ and $\reals^2\setminus K$ consists of a single component, then 
the level set flow of $K$ vanishes instantly.


\section{Smoothness of the level set flow}~\label{smoothness_section}

In this section we complete the proofs Theorem~\ref{main} and Theorem~\ref{jordan_curve} by showing that the boundary components of $K_t$ are smooth for all $t>0$.

\begin{thm} \label{boundary_smooth} Let $K\subset\reals^2$ be locally-connected, connected and compact, and let $\Gamma$ be a component of $\reals^2\setminus K$.  Then  $\partial\Gamma_t$ is a smooth closed curve for all $t>0$ such that $\Gamma_t\neq\emptyset$.
\end{thm} 

\begin{rem} The evolution of $\Gamma$ has been defined in the discussion preceding Proposition~\ref{LSF}.
\end{rem}

\proof  Using Remark~\ref{rem_approx}, let $\gamma_n$ be a sequence of smooth closed curves satisfying the conclusions of Proposition~\ref{approx_bound}.  In particular, $M_r(\gamma_n)<\widetilde{C}(r)$ for sufficiently large $n$.

Let $t>0$ be such that $\Gamma_t\neq\emptyset$, let $d=\diam(K)+1$ and let $r=r(t,d)$ be the scale chosen in Theorem~\ref{main_estimate}.  

{\bf Claim:}  For any $c>0$, there exists  $\nu>0$ such that for each $x\in\reals^2$ and $n$ sufficiently large,
each component of $(\gamma_n)_t\cap B_{\nu r}(x)$ is a $c$-Lipschitz graph over some line.  

This immediately implies that $\partial\Gamma_t$ is $C^1$ and Theorem~\ref{boundary_smooth} follows since $C^1$ curves instantly become smooth.

To prove the claim let $\gamma$ be the subcurve of $(\gamma_n)_t\cap B_{\nu
r}(x)$ that turns through the largest angle, i.e. $\gamma$ maximizes $|\int_\gamma\kappa ds|$.  The version of Corollary~\ref{curvature_estimate} mentioned in Remark~\ref{rem_approx} implies that after passing to a subsequence we may assume that $\int_{(\gamma_n)_t}
\kappa^2ds<C'$ for all $n$.  Then 
$$
\left(\int_\gamma\kappa ds\right)^2\leq \mathscr{L}(\gamma)\int_\gamma
\kappa^2ds<\nu C\widetilde{C}C',
$$
where the first inequality uses Cauchy-Scwartz, and the second uses Theorem~\ref{main_estimate} (2).  The claim follows by choosing $\nu<\left(2\arctan(c)\right)^2(C\widetilde{C}C')^{-1}$.\qed

\begin{cor} \label{boundary_main} Let $K\subset\reals^2$ be locally-connected, connected and compact.  Then $\partial K_t$ consists of finitely many disjoint smooth closed curves for each $t>0$.
\end{cor}

\proof Using Proposition~\ref{LSF}, Theorem~\ref{boundary_smooth} implies that $\partial K_t$ is a union of smooth closed curves, which are disjoint by the maximum principle, and the number of components is finite by Proposition~\ref{area}. \qed

We now have all the ingredients to complete the proofs of the main results concerning the level set flow:

\proof [Proofs of Theorems~\ref{main} and~\ref{jordan_curve}]  If $m(K)>0$ or $\reals^2\setminus K$ contains more than two components, then Proposition~\ref{area} implies that $m(K_t)$ is positive for small $t$.  Hence by Corollary~\ref{boundary_main}, the interior of $K$ is nonempty.  If the $m(K)=0$ and $\reals^2\setminus K$ consists of a single component, then simple area considerations show that $K_t$ vanishes instantly. 

The remaining case is when $\reals^2\setminus K$ consists of exactly two components.  According to Theorem~\ref{boundary_smooth}, $K_t$ lies between two smooth closed curves for small~$t$.  If $m(K)>0$ then by Proposition~\ref{area}, $m(K_t)>0$ for small $t$ and so the two boundary curves are distinct.  On the other hand, if $m(K)=0$ then $m(K_t)=0$ and the two curves coincide.\qed


\section{Backwards convergence}\label{backwards_convergence_section}

Since the level set flow of a measure zero Jordan curve $u=u_0$ is a smooth closed curve, any smooth parametrization of $u_t$, for some fixed $t>0$, extends to a solution of CSF, $u:\mathbb{S}^1\times(0,T)\to\reals^2$, whose image coincides with the level set flow.  A simple barrier argument shows that $u_t\subset\mathcal{N}_{\sqrt{2t}}(u)$, and it follows that $u_t$ Hausdorff converges to $u$ as $t\to 0$.  In this section we strengthen this convergence, showing that the level set flow converges as  a curve to the initial data.  More precisely:

\begin{thm}\label{backwards_convergence} Let $u:\mathbb{S}^1\to\reals^2$ be a measure zero Jordan curve and let $u_t$ be the level set flow of $u$.
For any $\epsilon>0$ there exists $T>0$ such that for any $0<t<T$ there exists a parametrization of $u_t$ such that $d_{\sup}(u,u_t)<\epsilon$.
\end{thm}

\begin{rem} For positive area Jordan curves the level set flow is the annular region between two smooth curves evolving by CSF.  The above statement is also true, and the same proof works in that case, if one replaces the level set flow with one of its boundary components.
\end{rem}

The proof of Theorem~\ref{backwards_convergence} proceeds by showing that if there were no such parametrization then the $r$-multiplicity would need to increase at some particular line, contradicting the following fact:

\begin{lem} \label{decreasing_mult}  Let $u$ be a measure zero Jordan curve.  Then $M_{r,\ell}(u_t)$ is a nonincreasing function of $t$.
\end{lem}

\proof This is clear when $t>0$ since then $u_t$ is smooth and the maximum principle applies.  Let $\{\gamma_n\}$ be a sequence of smooth closed curves converging uniformly to $u$.  Then $M_{r,\ell}(\gamma_n)\leq M_{r,\ell}(u)$ for large $n$.  If the $r$-multiplicity of the level set flow increases initially, then we obtain a contradiction to the fact that $M_{r,\ell}((\gamma_n)_t)$ is decreasing in $t$ since $(\gamma_n)_t$ Hausdorff converges to $u_t$.\qed

Let $\Omega$ be an annular neighbourhood of a Jordan curve $u\subset\reals^2$.  Suppose that $x_n, y_n\in u_{t_n}$ are two sequences of distinct points which limit onto the same point $x_0\in u$.  Let $[x_n,y_n]$ be the arc in $u_{t_n}$ which can be made into a null-homotopic loop in $\Omega$ by joining the endpoints of $[x_n,y_n]$ by a path contained in a small neighbourhood of $x_0$.  With this convention we have:

\begin{lem}\label{no_backtracks} Let $t_n\to 0$ be a sequence of positive times, and let $x_n, y_n\in u_{t_n}$ such that $\lim_{n\to\infty}x_n=\lim_{n\to\infty}y_n$.  Then 
$$
\lim_{n\to\infty}\diam([x_n,y_n])=0.
$$
\end{lem}

\proof If not, then after passing to a subsequence we may assume that 
$$
\diam([x_n,y_n])\geq\epsilon>0.   
$$
After passing to a further subsequence we may assume that there is a line~$\ell$ such that for each $n$, the 
arc $[x_n,y_n]$ intersects both components of $\reals^2\setminus\mathcal{N}_{r/2}(\ell)$, where $r=\frac{\epsilon}{2}$.  Moreover, we may assume that $\ell$ is chosen so that the components of $u$ defining $M_{r,\ell}(u)$ have distinct endpoints in $\partial\mathcal{N}_{r/2}(\ell)$.  This can be done since there are only finitely many lines in any given direction which do not have this property.

For each $\delta>0$, there exists smooth closed curves  $\phi_1,\phi_2:\mathbb{S}^1\to\reals^2$ lying in distinct components of $\reals^2\setminus u$ with $d_{\sup}(\phi_i,u)<\delta$ for $i=1,2$.  Let $\Omega$ be the region between $\phi_1$ and $\phi_2$, and  let $\{\Omega_i\}_{i\in\Lambda}$ be the components of $\Omega\cap\mathcal{N}_{r/2}(\ell)$ with the property that $\partial\Omega_i$ intersects both components of $\partial\mathcal{N}_{r/2}(\ell)$.  Then $|\Lambda|=M_{r,\ell}(u)$ for sufficiently small $\delta>0$ since the components defining $M_{r,\ell}(u)$ have distinct endpoints.

For each $t>0$ such that $u_t\subset\Omega$ there is an arc of $u_t$ that passes through each $\Omega_i$ since $u_t$ is homotopically nontrivial in $\Omega$.  Hence $M_{r,\ell}(u_t)\geq M_{r,\ell}(u)$ for such $t$.  But by assumption $[x_n,y_n]$ backtracks through $\Omega$ and hence must pass through some $\Omega_i$ at least twice.  Thus $M_{r,\ell}(u_t)>M_{r,\ell}(u)$ contradicting Lemma~\ref{decreasing_mult}. \qed

\proof[Proof of Theorem~\ref{backwards_convergence}] Fix $\epsilon>0$.  Let $\Omega$ be an annular neighbourhood of $u$ equipped with a (nonsmooth) foliation $\mathcal{F}$ such that 
\newline\indent(i) the endpoints of each leaf of $\mathcal{F}$ lie in distinct components of $\partial\Omega$,
\newline\indent(ii) each leaf of $\mathcal{F}$ intersects $u$ exactly once, and
\newline\indent(iii) the diameter of each arc is less than $\frac{\epsilon}{4}$.   
\newline Note that such a foliation exists by the strong Jordan Curve Theorem.

For each $t>0$ such that $u_t\subset\Omega$ there is a continuous map $\Phi_t:u_t\to u$ defined so that $x$ and $\Phi_t(x)$ belong to the same leaf of $\mathcal{F}$.  We say that $\Phi_t$ contains a backtrack of size $\epsilon$ if there exists an arc $[x,y]\subset u_t$ such that $\Phi_t(x)=\Phi_t(y$), $\diam(\Phi_t([x,y]))\geq\epsilon$ and the loop formed by taking $[x,y]$ and joining the endpoints along their common leaf in $\mathcal{F}$ is homotopically trivial in~$\Omega$. 

{\bf Claim:} If $\Phi_t$ contains no backtracks of size $\epsilon$, then there exists a parametrization of $u_t$ such that $d_{\sup}(u,u_t)<6\epsilon$.

Indeed, let $\{\alpha_i\}_{i\in\Delta}$ be a maximal ordered set of leaves of $\mathcal{F}$ such that $d(\alpha_i,\alpha_{i+1})\geq\epsilon$ for each $i$, where $d$ is the usual distance between closed sets, and let $A_i$ be the region in $\Omega$ between $\alpha_i$ and $\alpha_{i+1}$.   Note that maximality implies that $\diam(A_i)<3\epsilon$ for each $i\in\Delta$. 

Since we assume that $u_t$ contains no backtracks of size $\epsilon$ there is a subdivision of $u_t$ into $|\Delta|$ arcs $\{\gamma_i\}$ such that the endpoints of $\gamma_i$ lie in $\alpha_i$ and $\alpha_{i+1}$, and $\gamma_i$ is contained in $A_{i-1}\cup A_i\cup A_{i+1}$.  By identifying $\gamma_i$ with $u\cap A_i$ through a homeomorphism we obtain the desired parametrization.  This proves the claim.   

If the statement of the Theorem is not true for $6\epsilon$, then $\Phi_t$ contains a backtrack of size $\epsilon$ for arbitrarily small $t>0$.  By choosing a sequence $t_n\to 0$ so that the endpoints of the arcs causing such backtracks converge we obtain a contradiction to Lemma~\ref{no_backtracks}.\qed


\section{Finite length initial data}\label{cauchy_section}

In this section we show that the Cauchy problem for CSF has a unique solution whenever the initial data is a finite length Jordan curve.   

\begin{thm}\label{cauchy_problem} Let $J\subset\reals^2$ be the image of a finite length Jordan curve.  Then there exists $u\in\mathcal{C}^0(\mathbb{S}^1\times[0,T],\reals^2)\cap\mathcal{C}^\infty(\mathbb{S}^1\times(0,T],\reals^2)$ such that
\newline\indent (1) $u(\cdot,0)$ is a parametrization of $J$, and
\newline\indent (2) $u$ is a smooth solution to CSF on $\mathbb{S}^1\times (0,T)$. 

Moreover, the solution is unique up to reparametrization of $\mathbb{S}^1$.
 \end{thm}

The smoothness of the level set flow implies that the uniqueness part of Theorem~\ref{cauchy_problem} follows directly from the uniqueness of solutions to CSF for smooth closed curves.  As in the previous section let $u:\mathbb{S}\times (0,T)\to\reals^2$ be a smooth solution to CSF whose image coincides with the level set flow of $J$.  Since one can always choose an approximating sequence whose lengths converge to the length of the original curve the following fact is immediate:

\begin{lem}\label{length_limit} $\lim_{t\to 0}\mathscr{L}(u_t)=\mathscr{L}(J)$.
\end{lem}

The existence portion of Theorem~\ref{cauchy_problem} is proved by a series of lemmas.

\begin{lem}\label{length_piece} Let $x,y\in\mathbb{S}^1$and suppose $t_n\to 0$ is a sequence of times such that $\lim_{n\to\infty}u_{t_n}(x)=x_0$ and $\lim_{n\to\infty}u_{t_n}(y)=y_0$.  Then 
$$
\lim_{n\to\infty}\mathscr{L}(u_{t_n}[x,y])=\mathscr{L}([x_0,y_0]),
$$
where $[x_0,y_0]$ is the appropriate arc between $x_0$ and $y_0$ in $J$.
\end{lem}


\proof  Let $I_1$ and $I_2$ be the two arcs of $J$ between $x_0$ and $y_0$.  If $x_0=y_0$ then take $I_1=J$ and $I_2=x_0$.  Let $A_n=u_{t_n}[x,y]$ where $[x,y]$ is a choice of interval between $x$ and $y$ in $\mathbb{S}^1$, and let $B_n=u_{t_n}\setminus A_n$.  We may assume that $A_n$ is converging to $I_1$.  Then  $\liminf_{n}\mathscr{L}(A_n)\geq\mathscr{L}(I_1)$ and $\liminf_{n}\mathscr{L}(B_n)\geq\mathscr{L}(I_2)$.


But then Lemma~\ref{length_limit} implies that 
$$
\lim_{n\to\infty} \left(\mathscr{L}(A_n)+\mathscr{L}(B_n)\right)=\mathscr{L}(J),
$$
and hence that
\begin{align*}
\limsup_{n}\mathscr{L}(A_n)&=\mathscr{L}(J)-\liminf_{n}\mathscr{L}(B_n).
\\ &\leq\mathscr{L}(J)-\mathscr{L}(I_2)=\mathscr{L}(I_1),
\end{align*}
which proves the result.\qed

\begin{lem}\label{limit}  $\lim_{t\to 0} u_t(x)$ exists for each $x\in\mathbb{S}^1$.
\end{lem}

\smallskip

\proof We begin by showing that $\lim_{t\to 0}u_t(y)$ exists for some $y\in\mathbb{S}^1$.  Fix $T>0$ so that the $u_T$ exists, and for each $x\in\mathbb{S}^1$ define the path $\rho_x(t)=u_t(x)$ for $0<t<T$.  If $\lim_{t\to 0}\rho_x(t)$ does not exists then $\int_{\rho_x}|\kappa|dt=\infty$ since this quantity is the length of $\rho_x$.  This implies that
$$  
\int_{\rho_x}\kappa^2dt\geq \left(\int_{\rho_x}|\kappa|dt\right) -T=\infty.
$$

Now,  
$$
\mathscr{L}(u_t)-\mathscr{L}(u_T)=\iint\limits_{[t,T] u_t(\cdot)}\kappa^2 \diff s\diff t
$$
since the inner integral is the time derivative of length under CSF. Together with Lemma~\ref{length_limit} this implies that

$$
\mathscr{L}(J)\geq\iint\limits_{(0,T] u_t(\cdot)}\kappa^2 \diff s\diff t.
$$

The length of the tangent, $|\frac{\partial u}{\partial\theta}|$, is decreasing in $t$, since $u_t$ is evolving by CSF, and since there is a well-defined minimum on $u_T$ we have

$$
C=\inf_{(0,T]\times\mathbb{S}^1}\left\{\left|\frac{\partial u}{\partial\theta}(t,\theta)\right|\right\}>0.
$$

Thus
$$
\iint\limits_{(0,T] u_t(\cdot)}\kappa^2 \diff s\diff t=\iint\limits_{(0,T] [0,2\pi]}\kappa^2\left|\frac{\partial u}{\partial\theta}\right|\diff\theta\diff t> C\iint\limits_{(0,T] [0,2\pi]}\kappa^2\diff\theta\diff t,
$$
where the equality is a change of variables. 

By switching the order of integration we obtain
$$
\mathscr{L}(J)>\iint\limits_{[0,2\pi]\: \rho_x}\kappa^2\diff t\diff\theta=\infty,
$$
a contradiction.

Thus, there exists $y\in \mathbb{S}^1$ such that $\lim_{t\to 0}u_t(y)=y_0$ exists.  If $\lim_{t\to 0}u_t(x)$ does not exist for some $x\in\mathbb{S}^1$, then there exists a nontrivial subarc $[x_0,x_1]\subset J$ such that each point in $[x_0,x_1]$ is a limit point of the path $u_t(x)$, and $[x_0,x_1]$ may be chosen so that $y_0\notin [x_0,x_1]$.

Let $[y_0,x_0]$ and $[y_0,x_1]$ be subarcs of $J$ chosen so that $[y_0,x_0]\subset[y_0,x_1]$.  Then let $L_1=\mathscr{L}([y_0,x_0])$ and $L_2=\mathscr{L}([y_0,x_1])$, and note that $L_1<L_2$ since $[x_0,x_1]$ is nontrivial.  Then
$$
\limsup_{t\to 0}\mathscr{L}(u_t[y,x])\geq L_2
$$
since there exists a sequence of times $t_n$ such that $u_{t_n}[y,x]$ can be parametrized to converge uniformly to $[y,x_1]$. 

However, for any sequence of times $t_n\to 0$ such that $\lim_{n\to\infty}u_{t_n}(x)=x_0$, Lemma~\ref{length_piece} implies that 
$$\mathscr{L}(u_{t_n}[y,x])<\frac{L_1+L_2}{2}
$$
for large $n$.  However, since $u_t$ evolves by CSF, the length of such a segment is decreasing, and therefore  
$$
\mathscr{L}(u_t[y,x])<\frac{L_1+L_2}{2}<L_2
$$ 
for all $t>0$, a contradiction.\qed

Lemma~\ref{limit} implies that $u$ can be extended to a function on $\mathbb{S}^1\times[0,T]$ by setting $u(x,0)=u_0(x)=\lim_{t\to 0}u_t(x)$.  Theorem~\ref{cauchy_problem} then follows from next two results:

\begin{lem} $u_0$ is continuous.
\end{lem}

\proof Suppose $x_n\to x$ in $\mathbb{S}^1$.  Let $I_n$ be an interval in $\mathbb{S}^1$ between $x_n$ and $x$, chosen so that $\diam(I_n)\to 0$ as $n\to\infty$, and let $J_n$ be the interval in $J$ between $u_0(x_n)$ and $u_0(x)$ which contains $u_0(y)$ for each $y\in I_n$. 

Given $\epsilon>0$, Lemma~\ref{length_limit} implies that there exists $T>0$ such that  $\mathscr{L}(J)-\mathscr{L}(u_t)<\epsilon$ for all $0<t<T$.  And since $u_T$ is a smooth closed curve and $\diam(I_n)\to 0$ we can choose $n_0\in\naturals$ such that $\mathscr{L}(u_T(I_n))<\epsilon$ for all $n\geq n_0$.  Thus 
 $$
 \mathscr{L}(u_T(\mathbb{S}^1\setminus I_n))>\mathscr{L}(J)-2\epsilon,
$$
and the same inequality holds for $0<t\leq T$ since $u_t$ is evolving by CSF.

Applying Lemma~\ref{length_piece} to $u_t(\mathbb{S}^1\setminus I_n)$ we obtain that
$$
\mathscr{L}(J\setminus J_n)\geq\mathscr{L}(J)-2\epsilon.
$$  

Thus
$$
\mathscr{L}(J_n)\leq2\epsilon,
$$ 
and hence $d(u_0(x_n),u_0(x))\leq 2\epsilon$ for all $n\geq n_0$.\qed

\begin{lem} $u_0$ is injective.
\end{lem} 

\proof  Lemma~\ref{backwards_convergence} implies that if $u_0(x)=u_0(y)$, then for one of the two intervals, say $[x,y]$, between $x$ and $y$ in $\mathbb{S}^1$ we have $ u_0(z)=u_0(x)$ for all $z\in[x,y]$.  Lemma~\ref{length_piece} then implies that $\lim_{t\to 0}\mathscr{L}(u_t[x,y])=0$.  But $\mathscr{L}(u_t[x,y])$ is decreasing in $t$ and positive for $t>0$, and this is impossible.\qed



\addcontentsline{toc}{section}{References}

\end{document}